\documentclass[13pt]{article}

\usepackage{amsmath}
\usepackage{amsbsy}
\usepackage{amsthm}
\usepackage{amssymb}
\usepackage{amsfonts}
\usepackage{dsfont}
\usepackage{color}

\newtheorem{theorem}{Theorem}

\newtheorem{lemma}[theorem]{Lemma}

\def\E{{\mathds E}}
\def\I{{\mathds 1}}

\def\P{{\mathds P}}
\def\R{{\mathds R}}

\def \beq{\begin{equation}}
\def \eq{\end{equation}}
\def \proba{\xrightarrow[n]{(proba.)}}
\def \eref#1{(\ref{#1})}

\def\z{{\zeta}}

\newcommand{\conv}{\textrm{conv}}

\def \Proof{{\bf Proof. }}

\begin{document}

\title{Many empty triangles have a common edge}

\author{Imre B\'ar\'any\thanks{R\'enyi Institute of Mathematics, Hungarian Academy of Sciences PO Box 127, 1364 Budapest  (Hungary), and
 Department of Mathematics, University College London, Gower Street,  London WC1E 6BT (England). Partially supported by ERC Advanced Research Grant no 267165 (DISCONV),
and by Hungarian National Research Grants No K84767 and NK78439}
 \and  Jean-Fran\c{c}ois Marckert \thanks{LaBRI, CNRS, Universit\'e Bordeaux,
351 cours de la Lib\'eration, 33405 Talence cedex (France), partially supported by ANR blanc PRESAGE (ANR-11-BS02-003)} \and  Matthias Reitzner,\thanks{Institut f\"ur Mathematik, Universit\"at Osnabr\"uck, 49069 Osnabr\"uck (Germany)}}

\maketitle
\begin{abstract} Given a finite point set $X$ in the plane, the degree of a pair $\{x,y\} \subset X$ is the number
of {\sl empty triangles} $t=\conv\{x,y,z\}$, where empty means $t\cap X=\{x,y,z\}$. Define $\deg X$ as the maximal
degree of a pair in $X$. Our main result is that if $X$ is a random sample of $n$ independent and uniform points
from a fixed convex body, then $\deg X \ge cn/\ln n$ in expectation.

\end{abstract}

\subsubsection*{Acknowledgements} The authors thank the organizers of the workshop PRESAGE (Cluny, France 2012) for arranging such an enjoyable and fruitful conference, where this work got started.

\section{Introduction}
Let $X$ be a finite set of points in $\mathbb{R}^2$ in general position meaning that no 3 points of $X$ are collinear. For $k\geq 2$,  let ${X \choose k}$ be the family of all $k$-element subsets of $X$.  The triangle $t:=\conv \{x_1,x_2,x_3\}$ where $\{x_1,x_2,x_3\} \in {X \choose 3}$ is said to be empty (in $X$) if $t \cap X= \{x_1,x_2,x_3\}$.

The degree of a pair $\{x,y\} \in {X \choose 2}$ is the number of $z \in X$ such that $\{x,y,z\} \in {X \choose 3}$ determines
an empty triangle. The degree of $\{x,y\}$ will be denoted by $\deg (x,y)$ or $\deg (x,y;X)$. We set
\[
\deg X= \max \left\{ \deg (x,y): \{x,y\} \in {X \choose 2} \right\}.
\]
Setting $|X|=n$ we clearly have $\deg X \le n-2$. The following conjectures was raised by the first author, and appeared first in a paper of Paul Erd\H os \cite{er} in 1992 and repeated in \cite{bk}.

\medskip
{\bf Conjecture 1.} $\deg X $ goes to infinity as $n\to \infty$.
\medskip

Very little is known about the validity of the conjecture. Namely, it is shown in \cite{bk} that $\deg X\ge 10$ for large enough $n$. The construction in \cite{bv} gives a set $X \in \mathbb{R}^2$ in general position with $\deg X =4\sqrt n(1+o(1))$, as one can check easily. 

In this paper we give a lower bound on $\deg X$ in a special case, namely, when $X=\xi_n$ is a set of $n$ independent, random points chosen uniformly from a fixed compact convex set $C \subset \R^2$ that has nonempty interior. Then $\xi_n$ is in general position with probability one, so $\deg \xi_n$ is a well-defined random variable. Our main result shows that the expectation of $\deg \xi_n$ is quite close to $n$.

\begin{theorem}\label{th:Eet} There is a universal constant $c>0$ such that
$$ \E (\deg \xi_n)  \geq \frac {cn}{\ln n} . $$
\end{theorem}

The following theorem is not a consequence of the previous one, even though it seems weaker.
\begin{theorem}\label{th:Eet2}
\[\deg \xi_n \proba +\infty.\]
\end{theorem}

Although Theorem  \ref{th:Eet2} can be proved by a modification of the argument used for Theorem \ref{th:Eet} we provide an independent proof relying on a ``local argument'' that can be used quite generally. This local argument is the following. Take a grid with mesh $1/\sqrt{n}$ on the plane. This grid defines certain squares, and each square $Q$ which is totally included in the underlying convex body $C$ contains a ${\sf Binomial}(n,1/n)$ number $N_Q$ of points of $\xi_n$. When $n$ is large, $N_Q$ can be approximated by a Poisson random variable with parameter 1. Even if these variables $N_Q$ are not independent, some results can be transferred from the case where they are (this is in substance Lemma \ref{lem:transfer}). Then, using a collection of ''$n$ independent squares'', each of them containing Poisson(1) number of random points, it is not difficult to see that one of these squares will contain a set $Y$ of $k$ points satisfying $\deg Y=k-2$, with probability  going to 1, since the probability that one of them satisfies this condition is positive (see details in the proof of Theorem \ref{th:Eet2}, below).

This local argument gives the following more general statement concerning "order types" \cite{gp}. Two finite sets $A,B \subset \R^2$ are of {\sl same type} if both are in general position and there is a one-to-one correspondence between $A$ and $B$, say $a_i \leftrightarrow b_i$ ($i=1,\dots,k$) where $A=\{a_1,\dots,a_k\}$ and $B=\{b_1,\dots,b_k\}$, such that the orientations of the triangles $(a_h,a_i,a_j)$ and $(b_h,b_i,b_j)$ are the same for all $1\le h<i<j\le k$. Being of the same type is an equivalence relations and equivalence classes are called {\sl order types}.

\begin{theorem} Let $P$ be a fixed order type. Then, as $n \to \infty$,
\[ \P( \xi_n \cap Q \mbox{ \rm is of type } P \mbox{ \rm for some square }Q) \to 1.
\]
\end{theorem}

The proof goes the same way as that of Theorem \ref{th:Eet2} and is therefore omitted.

\smallskip
We close this section by stating (or rather repeating) another conjecture from \cite{bk}. Given a finite $X \in \R^2$ in general position, let $f(X)$ denote the number of empty triangles in $X$. Set
\[
f(n)=\min \{f(X): X \subset \R^2 \mbox{ is in general position and } |X|=n\}.
\]
It is known that
\[n^2-5n \le f(n) \le 1.6195\dots n^2
\]
where the lower bound is from \cite{gar} and the upper one from \cite{bv}. We think that the lower bound is larger than $n^2$. We state the following

\smallskip
{\bf Conjecture 2.} For large enough n, $f(n) > 1.01n^2$.

\smallskip
The same conjecture appears also in \cite{bk}. Note that, when $\xi_n$ is the random uniform sample from $C$, the expectation of $f(\xi_n)$ is of order $2n^2$, see Valtr~\cite{va}.

\section{Proof of Theorem \ref{th:Eet} }

The proof relies of a first moment argument that we state for arbitrary $X$ first. For any $T>0$, any set $X$, we have
\beq\label{eq:first}
\sum\limits_{ \{x,y\} \in {X \choose 2}} \I( \| x-y \| \leq T )\, \deg (x,y;X) \leq  N_T (X) \deg X
\eq
where
\[N_T (X)=\sum_{ \{x,y\} \in {X \choose 2}} \I( \| x-y \| \leq T )\]
is the number of pairs $\{x,y\} \in {X \choose 2}$ with $\| x-y \| \leq T $.  From \eref{eq:first}, we get the following formula which will play a central role here:
\beq\label{eq:first2}
\deg X \geq \frac{1}{N_T(X)} \sum_{ \{x,y\} \in {X \choose 2}} \I( \| x-y \| \leq T ) \deg (x,y).
\eq
A comment must be added here : when $N_T(X)=0$, the division by $N_T (X)$ is not valid in \eref{eq:first2}. Notice however that in this case $\sum_{ \{x,y\} \in {X \choose 2} }\I( \| x-y \| \leq T ) \deg (x,y;X)$ is also 0, and then the right hand side in \eref{eq:first2} has the form 0/0. Thus in this case we consider this ratio to be 0, and then \eref{eq:first2} is indeed valid (the problem disappears in \eref{eq:1KE}).

We need some geometric preparations. Clearly $\deg X$ is invariant under non-degenerate affine transformation (and so is $f(X)$). Also, $\xi_n$ is invariant (or rather equivariant) under such a transformation. So we can apply an arbitrary (non-degenerate) transformation to $C$, and $\E \deg \xi_n$ will not change. By the L\"owner half of the L\"owner-John theorem (see for instance \cite{jo}, or \cite{ba} for a more modern treatment), there is a pair of concentric ellipses $E_1,E_2$ so that $E_2$ is a blown up copy of $E_1$ by a factor of 2 with $E_1 \subset C \subset E_2$. First we apply the affine transformation so that the area of $C$ becomes equal to one. This is convenient since then the Lebesgue measure coincides with the probability measure defining $\xi_n$. Second we apply an area preserving affine transformation that carries $E_1$  resp. $E_2$ to $rD$ and $2rD$ where $D$ the Euclidean unit disk, centered at the origin. It is easy to see that one can take $r=27^{-1/4}$. From now on we assume that $C$ is in this position, and $rD \subset C \subset 2rD$.

We can put now $X=\xi_n$ in \eref{eq:first2} and take the expectation on both sides. This yields
\beq \label{eq:first4}
\E(\deg \xi_n) \geq  \E\left(\frac{1}{N_T(\xi_n )} \sum_{ \{x,y\} \in {\xi_n \choose 2}} \I( \| x-y \| \leq T ) \deg (x,y;\xi_n)\right),
\eq
and further, for any $K>0$,
\beq \label{eq:1KE}
\E \left(\deg \xi_n \right)  \geq
\frac{1}{K} \ \E \left[ \sum_{ x,y \in {\xi_n \choose 2}} \I( \| x-y \| \leq T ) \deg (x,y; \xi_n) \I(N_T(\xi_n) \leq  K)
\right].
\eq
The idea now is to somehow optimise in $T=T_n$ and in $K=K_n$. The intuition here is that in the random case
$\deg \xi_n$ will be reached for a pair $\{ x,y\}$ with $x,y$ very close. The best we can do here is to take $T=1/n$ and $K=c \ln n$ for some $c>0$ to be fixed later.

{\bf Remark.} The closest pair in ${\xi_n \choose 2}$ will have distance approximately $n^{-1}$ and we expect the closest pair of points in $\xi_n$ to be very likely to give the maximal, or close to the maximal degree. In fact, for $T:=T_n\to 0$, it is a well known fact (see, e.g., the recent paper by Reitzner, Schulte, and Thaele~\cite{rst}) that
$$ \E \left(N_T (\xi_n)\right)  \approx \frac {\pi}2 n^2  T^2  \left( 1+ o(1) \right) $$
implying that the expected number of pairs of distance $n^{- 1}$ is positive.

Before choosing these special values for $T$ and $K$, we go on from \eref{eq:1KE} by starting the computation of the right hand side, to be denoted by ${\sf RHS}$ from now on. By conditioning successively on each pair $\{x_i,x_j\}$ of ${\xi_n \choose 2}$, we have
\begin{equation} \label{eq:RHS2}
{\sf RHS}=\frac{{n \choose 2}}{K}
\int \limits_{C^2} \I(  \| x-y \| \leq T ) \E \left[ \deg (x,y;\xi_{n-2}) \I(N_T(\xi_{n-2}\cup\{ x,y\}) \leq K) \right] \, dy dx,
\end{equation}
here $\deg (x,y;\xi_{n-2})$ denotes the number of empty triangles with base $\{x,y\}$ in
 $\z:= \xi_{n-2}\cup\{x,y\}$.

Since
$$\E (X \I(A) ) = \E( X) - \E( X \I(A^c)) $$
we have the following bound on the factor in \eref{eq:RHS2}
\begin{eqnarray}
\E \left[ \deg (x,y;\xi_{n-2}) \I(N_T(\z) \leq K) \right]
&=& \nonumber
\E \left[ \deg (x,y;\xi_{n-2}) \right]  - \E \left[ \deg (x,y;\xi_{n-2}) \I(N_T(\z) > K) \right]
\\ & \geq & \label{eq:E-E}
\E \left[ \deg (x,y;\xi_{n-2}) \right] - \E \left[ n \I(N_T(\z) \geq K)  \right].
\end{eqnarray}

Choose now $\rho=r/2$. Then the disk $\rho D$ is contained in $rD \subset C$ and every point of $\rho D$ is farther than $\rho=r/2$ from the boundary of $C$.

\begin{lemma} \label{lem:ad} If $\| x-y\| \leq T = T_n:=n^{-1}$, and $x,y \in \rho D$, then $ \E\left[ \deg (x,y;\xi_{n-2}) \right] \geq   \rho n (1-e^{-\rho /2}). $
\end{lemma}

We need one more ingredient, actually a crucial one, which is a special case of Theorem 5.3 from \cite{rst}.
\begin{lemma}\label{lem:crux} Assume $\alpha >0$ and  $T=\alpha n^{-1}$. Then there is $c=c(\alpha)>0$ such that
\[
 \P( N_{T} (\xi_n) \geq 144 \ln n) \leq c n^{-3}.
 \]
\end{lemma}

Assuming for a moment that Lemma~\ref{lem:ad} has been proved, we finish the proof of Theorem \ref{th:Eet}. Note first that
\begin{eqnarray}
\E \left[ n \I(N_T(\z) \geq K)  \right] = n\P \left[ N_T(\z) \geq K \right].
\end{eqnarray}
Further, let $D(x,T)$ and $D(y,T)$ denote the disks centered at $x$ and $y$ with radius $T$. Then
\[
N_T(\z)\le N_T(\xi_{n-2})+|\xi_{n-2}\cap D(x,T)|+|\xi_{n-2}\cap D(y,T)|+1 \le   N_T(\xi_{n-2})+2 N_{2T}(\xi_{n-2})+1 \le 3 N_{2T}(\xi_n)+1
\]
where we used the fact that if two points lie in $D(x,T)$, then their distance is at most $2T$. This implies, with $K=K_n=3\cdot 145\ln n$, that
\[
\P \left[ N_T(\z) \geq K_n \right]\le \P \left[ 3N_{2T}(\xi_n)+1 \geq K_n \right] \le \P \left[ N_{2T}(\xi_n) \geq 144\ln n \right].
\]
We set $T=T_n=1/n$. Then, using Lemma~\ref{lem:crux} with $\alpha=2$, we have with a suitable $c>0$
\[ \E [ n \I(N_{T_n}(\xi_{n-2}) \geq K_n) ] 
\le n\P [ N_{2T_n}(\xi_n) \geq 144 \ln n ]\le cn^{-2}.
\]

Hence we see that the second term in \eref{eq:E-E} is negligible compared to the first one, whose value is more than
$n\rho (1-e^{-\rho /2})> 0.086n$, by Lemma \ref{lem:ad}, for $n$ large enough.

Plugging what is known into \eref{eq:E-E}, we obtain
\[{\sf RHS}\geq  \frac{n(n-1)}{2K_n}  \left(n\rho (1-e^{-\rho /2})-o(1)\right) \int_{(\rho D)^2} \I(  \| x-y \| \leq T_n ) dx dy.\]
The value of the integral being larger than $c'/n^2$ for some $c'$ (in fact, it is equivalent to $\pi \rho^2 T^2_n$), since  $K_n$ is constant times $\ln n$, this ends the proof of Theorem \ref{th:Eet}. \hfill$\Box$

\medskip
Lemma \ref{lem:ad} remains to be proved.

\smallskip
\noindent{\bf Proof of Lemma \ref{lem:ad}.} Fix $x$ and $y$ in $\rho D$ with $\|x-y\|\le 1/n$. We have
\begin{eqnarray*}
\E \left[ \deg(x,y;\xi_{n-2}) \right]
&=&
\E \left[
\sum_{z \in \xi_{n-2}} \I( \xi_{n-2} \cap \conv \{x,y,z\}=\{x,y,z\} )
\right]
\\ & = &
(n-2)\P (  \xi_{n-3} \cap \conv \{x,y,U\}=\{x,y,U\} )
\end{eqnarray*}
where $U$ is a random variable uniform in $C$ independent from $\xi_{n-2}$. This gives
\begin{eqnarray*}
\E \left[ \deg (x,y;\xi_{n-2}) \right]   & = &
(n-2) \int\limits_{C} (1- A(x,y,u))^{n-3} du.
\end{eqnarray*}
where $A(x,y,u)$ is the area of the triangle $x,y,u$. (This is where the condition that the area of $C$ is one is convenient.)  Let $Q(x,y)$ be the square of side length $\rho$, centered at $(x+y)/2 \in \rho D$, with one side parallel with the vector $x-y$. Instead of integrating (with respect to $u$) on $C$, we integrate only on the square $Q(x,y)$. This gives the lower bound:
\begin{eqnarray*}
 \int\limits_{C} (1- A(x,y,u))^{n-3} du
& \geq & \int\limits_{Q(x,y)} (1- A(x,y,u))^{n-3} du  \\
&\geq &\int\limits_{[0,\rho]} \int\limits_{[0,\rho]} (1- \frac 12 \|y-x\| z_2)^{n-3} dz_2dz_1
\\ &=&
\rho \frac 2{\|x-y\|}  \int\limits_{[0, \|x-y\| \rho /2]} (1- t)^{n-3} dt
\\ &=&
\frac {2\rho}{\|x-y\|}  \frac 1{n-2} \ [ 1- (1- \|x-y\|\rho /2)^{n-2}  ]
\\ & \geq &  \frac {2\rho n}{n-2} \ [ 1- (1- \rho/(2n))^{n-2}  ]\\
&\geq &   \frac {2\rho n}{n-2} [1- e^{- (n-2)\rho /(2n)} ]  \geq  \frac {\rho n}{n-2} (1-e^{-\rho/2})
\end{eqnarray*}
for $n $ large enough. The last inequality is a consequence of $[1- e^{- (n-2)\rho /(2n)} ]\to 1-e^{-\rho/2}$.
\hfill $\Box$

\section{Proof of Theorem \ref{th:Eet2} }

Consider the convex body $C$ and $\xi_n$ as defined at the beginning of the paper (we don't need the disks $rD,2rD$ now). Add a grid with mesh $1/\sqrt{n}$ on the plane. The squares $(Q_i,i \in I)$ hence obtained have area $1/n$. For $n$ large enough (that we choose even, for convenience),  choose $n/2$ squares $Q_1,\dots,Q_{n/2}$ totally included in $C$, and consider $C^\star=C\setminus \cup_{i=1}^{n/2}Q_i$ be the remaining part of $C$.

For $1\leq i \leq n/2$, let $\xi_n^{(i)}=\xi_n\cap Q_i$, and  set $\xi_n^\star=\xi_n\cap C^\star$ the ``composition'' of the square of interests and the other ones.
Further, let $N_i=|\xi_n^{(i)}|$, and let $N^\star=|\xi_n\cap C^\star|$.\par
The family $(N_1,\dots,N_{n/2},N^\star)$ has multinomial distribution ${\sf Mult}(1/n,\dots,1/n,1/2)$. As such it is distributed as a collection of independent Poisson random variables $(P_1,\dots,P_{n/2},P^\star)$ conditioned by $P^\star+\sum_{i=1}^{n/2}P_i=n$, where the $P_i$'s are {\sf Poisson}(1) distributed and $P^\star$ is {\sf Poisson}($n/2$).

We first claim that rare events for Poisson random variables $(P_1,\dots,P_{n/2})$ are also rare events for the first $n/2$ marginals $(P_1,\dots,P_{n/2},P^\star)$ conditioned by $P^\star+\sum_{i=1}^{n/2}P_i=n$.
\begin{lemma}\label{lem:transfer}There exists a constant $b>0$, such that for large enough $n$ and for any measurable set $A\in \mathbb{R}^{n/2}$,
\beq
 \P\left((P_i,i=1,\dots,n/2) \in A \Big| P^\star+\sum_{i=1}^{n/2}P_i=n\right) \leq b\,\P((P_i,i=1,\dots,n/2)\in A).
\eq
\end{lemma}
\Proof  Take any $m_1,\dots,m_{n/2}\in \mathbb{N}$ and write
\[\P\left(P_i=m_i,i=1,\dots,n/2 \Big| P^\star+\sum_{i=1}^{n/2}P_i=n\right)=\P\left(P_i=m_i,i=1,\dots,n/2\right)\frac{\P(P^\star=n-\sum_{i=1}^{n/2} m_i)}{\P\left( P^\star+\sum_{i=1}^{n/2}P_i=n\right)}.\]
By Stirling, $\P\left( P^\star+\sum_{i=1}^{n/2}P_i=n\right)\sim 1/\sqrt{2\pi n}$, and at the numerator, it is easily checked that $\P(P^\star=n-\sum_{i=1}^n m_i)\leq \max_k \P(P^\star=k)=O(1/\sqrt{n})$, which ends the proof. \hfill $\Box$ \medskip

We now consider the same collection of $M=n/2$ squares $Q_1,\dots,Q_M$ as above, but in each of them we place
independently of the other ones, a Poisson point process $\zeta^{(j)}$ with intensity 1.
Here $P(j):=|\zeta^{(j)}|\sim {\sf Poisson}(1)$, and conditionally on $P(j)=k$, the $k$ points in $\zeta^{(j)}$  are taken in $Q_j$ according to the uniform distribution.

A simple application of the Lemma and of the fact that for both models the points in the squares are i.i.d. uniform (and independent of the other squares) given their numbers, we get that for any large enough $n$ and for any measurable set $A$
\[\P( (\xi_n^{(i)},i=1,...,n/2)\in A) \leq b\,  \P( (\zeta^{(i)},i=1,\dots,n/2)\in A),\]
with the same constant $b$ as in the Lemma.

 In the Poisson model we have more independence, and it is then more easy. The number of squares $Q_j$ such that $P(j)=k$ is ${\sf Binomial}(M,e^{-1}/k!)$ distributed. Let $\xi_k$ be a random, independent sample of $k$ points chosen uniformly from the unit square, and set
$p_k=\P(\deg \xi_k=k-2)$. Clearly $p_k>0$: for instance if the $k$ points are in a convex position (which happens with positive probability), then all triangles formed with 3 different points are empty. \par
It follows that the number of squares $Q_j$ such that $P(j)=k$ and $\deg (\zeta^{(j)})=k-2$ has law ${\sf Binomial}(M,(e^{-1}/k!) p_k)$. When $M$ goes to $+\infty$, at least one of these squares  satisfies $P(j)=k$ and $\deg(\zeta^{(j)})=k-2$ with probability going to 1.
This implies that  for any $L$, letting
\[B_L:=\{\forall l\leq L, \exists j : P(j)=l, \deg (\zeta^{(j)})=l-2\}\]
we have
\[\P(B_L)\to 1.\]
The probability of the complementary event $A=\complement B_L$ goes to 0. Lemma \ref{lem:transfer} allows then to see that
\[ \P\left((P_i,i=1,\dots,n/2) \in A \Big| P^\star+\sum_{i=1}^{n/2}P_i=n\right)\to 0\]
which concludes the proof of Theorem \ref{th:Eet2}.\hfill $\Box$

\end{document}